\documentclass[11pt]{amsart}
\usepackage{palatino,hhline, bbm}
\usepackage{amsthm, amsmath}
\usepackage{amssymb} 
\usepackage{amscd}
\usepackage{mathrsfs}
\renewcommand{\mathcal}{\mathscr}
\usepackage[hmargin=2.8cm, vmargin=2.3cm]{geometry}
\usepackage[breaklinks,bookmarksopen,bookmarksnumbered]{hyperref}
\usepackage[all]{xy}
 \usepackage[applemac]{inputenc}

\theoremstyle{plain}
\newtheorem{thm}{Theorem}
\newtheorem{lem}{Lemma}
\newtheorem{prop}{Proposition}
\newtheorem{cor}{Corollary}
\newtheorem{defi}{Definition}
\theoremstyle{remark}

\newcommand\pr{\noindent\textit{Proof} : }

\newcommand\rond{\kern 1pt{\scriptstyle\circ}\kern 1pt}
\newcommand\pt{\scriptscriptstyle \otimes }

\newcommand\Aut{\operatorname{Aut}}

\newcommand\im{\operatorname{Im}}
\newcommand\Ker{\operatorname{Ker}}
\newcommand\Coker{\operatorname{Coker}}
\newcommand\Pic{\operatorname{Pic}}

\newcommand\Tr{\operatorname{Tr}}
\newcommand\tor{\operatorname{Tors}}
\newcommand\Spec{\operatorname{Spec}}
\newcommand\Sing{\operatorname{Sing}}

\newcommand\pp{p.p.a.v.\ }

\newcommand\Z{\mathbb{Z}}
\newcommand\Q{\mathbb{Q}}
\newcommand\R{\mathbb{R}}
\newcommand\C{\mathbb{C}}
\renewcommand\P{\mathbb{P}}

\newcommand\F{\mathbb{F}}
\newcommand\G{\mathbb{G}}
\renewcommand\O{\mathcal{O}}

\def\qfl#1{\buildrel {#1}\over {\longrightarrow}}
\newcommand\mono{\lhook\joinrel\mathrel{\longrightarrow}}
\newcommand\iso{\vbox{\hbox to .8cm{\hfill{$\scriptstyle\sim$}\hfill}
\nointerlineskip\hbox to .8cm{{\hfill$\longrightarrow $\hfill}} }}
\newcommand\bir{\vbox{\hbox to .8cm{\hfill{$\scriptstyle\sim$}\hfill}
\nointerlineskip\hbox to .8cm{{\hfill$\dasharrow $\hfill}} }}
\def\sdir_#1^#2{\mathrel{\mathop{\kern0pt\oplus}\limits_{#1}^{#2}}}

\begin{document}
\title[The L\"uroth problem]{The L\"uroth problem}
%\subtitle{Preliminary version}
\author[Arnaud Beauville]{Arnaud Beauville}
\address{Laboratoire J.-A. Dieudonn\'e\\
UMR 7351 du CNRS\\
Universit\'e de Nice\\
Parc Valrose\\
F-06108 Nice cedex 2, France}
\email{arnaud.beauville@unice.fr}
 
%\date{\today}
\begin{abstract}
\vskip0.3cm The L\"uroth problem asks whether every unirational variety is rational. After a historical survey, we describe  the methods developed in the 70's to get a negative answer, and give some easy examples. Then we discuss a new method introduced last year by C. Voisin.
\end{abstract}

\maketitle

\section{Some history}
\subsection{Curves and surfaces}
In 1876 appears a three pages note by J. L\"uroth \cite{L}, where he proves that if a complex algebraic  curve $C$ can be parametrized by rational functions, one can find another parametrization which is generically one-to-one. In geometric language, if we have a dominant rational map $f:\P^1\dasharrow C$, then $C$ is a rational curve.

By now this is a standard exercise : we can assume that $C$ is smooth projective, then $f$ is a morphism, which induces an injective homomorphism $f^*: H^0(C,\Omega ^1_C)\rightarrow \allowbreak H^0(\P^1,\Omega ^1_{\P^1})=0$. Thus $C$ has no nontrivial holomorphic 1-form, hence has genus $0$, and this implies $C\cong \P^1$.

Actually L\"uroth does not mention at all  Riemann surfaces, but uses instead  an ingenious and somewhat sophisticated algebraic argument. I must say that I find somewhat  surprising  that he did not consider applying Riemann's theory, which had appeared 20 years before. 

Anyhow, clearly L\"uroth's paper had  an important impact.  
When  Castelnuovo and Enriques develop the theory of algebraic surfaces, in the last decade of the 19th century, one of the first questions they attack is whether the analogous statement holds for surfaces. Suppose we have a smooth projective surface $S$ (over $\C$) and a dominant rational map $f:\P^2\dasharrow S$. As in the curve case, this implies $H^0(S,\Omega ^1_S)=H^0(S,K_S)=0$ (note that $f$ is well-defined outside a finite subset). At first Castel\-nuovo hoped that this vanishing would be sufficient to characterize rational surfaces, but Enriques suggested a counter-example,  now known as the Enriques surface. Then Castelnuovo found the right condition, namely $H^0(S,\Omega ^1_S)=H^0(S,K_S^2)=0$; this is satisfied by our surface $S$, and Castelnuovo proves that it implies that $S$ is rational. After more than one century, even if the proof has been somewhat simplified, this is still a highly nontrivial result.

%\smallskip	
\subsection{Attempts in dimension 3}
At this point it becomes very natural to ask  what happens in higher dimension. Let us first
recall the basic definitions (see \S 3 for a more elaborate discussion): a complex variety $X$ of dimension $n$ is \emph{unirational} if there is a dominant rational map $\P^n\dasharrow X$; it is \emph{rational} if there is a birational such map. The \emph{L\"uroth problem} asks whether every unirational veriety is rational.

In 1912, Enriques proposed a counter-example in dimension 3 \cite{E}, namely a smooth complete intersection of a quadric and a cubic in $\P^5$ -- we will use the notation $V_{2,3}$ for such a complete intersection. Actually what Enriques does in this two pages paper is to prove the unirationality of  $V_{2,3}$, in a clever (and correct) way; for the non-rationality he refers to a 1908 paper by Fano \cite{F1}. 

In the course of his  thorough study of what we now call Fano manifolds, Fano made various attempts to prove that some of them are not rational \cite{F2,F4}. Unfortunately the birational geometry of threefolds is considerably more complicated than that of surfaces; while the intuitive methods of the Italian geometers were sufficient to handle surfaces, they could not treat adequately higher-dimensional manifolds.  None of Fano attempted proofs is acceptable by  modern standards. 

A detailed criticism of these attempts can be found in the book \cite{Roth}. It is amusing that after concluding that none of them can be considered as correct, Roth goes on and proposes a new counter-example, which is not simply connected and therefore not rational  (the fundamental group is a birational invariant). Alas, a few years later Serre (motivated in part by Roth's claim) proved that a unirational variety is simply connected \cite{Serre}.

\subsection{The modern era}\label{mod}
Finally, in 1971-72, three different (indisputable) counter-examples appeared. We will discuss at length these results in the rest of the paper;
let us indicate briefly here the authors, their examples and the methods they use to prove non-rationality :

\medskip	
\begin{center}\setlength{\tabcolsep}{6pt}

{\renewcommand{\arraystretch}{1.5}\begin{tabular}{|c|c|c|}\hline 
 Authors & Example & Method  \\
\hline 
Clemens-Griffiths & $ V_3\subset \P^4 $ & $ JV $ \\ 
\hline 
Iskovskikh-Manin & some $ V_4\subset \P^4 $ & $ \mathrm{Bir}(V) $ \\
\hline 
Artin-Mumford & specific & $ \mathrm{Tors}\ H^3(V,\Z) $ \\
\hline
\end{tabular}}
\end{center}
\medskip	

More precisely :

\smallskip	
$\bullet$ Clemens-Griffiths \cite{CG} proved the longstanding conjecture that a smooth cubic threefold $V_3\subset \P^4$ is not rational   -- it had long been known that it is unirational. They showed that the \emph{intermediate Jacobian}  of $V_3$ is not a Jacobian (\emph{Clemens-Griffiths criterion}, see Theorem \ref{CG} below).

\smallskip	
$\bullet$ Iskovskikh-Manin \cite{IM} proved that any smooth quartic threefold $V_4\subset \P^4$ is not rational. Some unirational quartic threefolds had been constructed by B. Segre \cite{Segre}, so these provide counter-examples to the L\"uroth problem. They showed that the group of birational automorphisms of $V_4$ is finite, while the corresponding group for $\P^3$ is huge. 

\smallskip	
$\bullet$ Artin-Mumford \cite{AM} proved that a particular double covering $X$ of $\P^3$, branched along a quartic surface in $\P^3$ with 10 nodes, is unirational but not rational. They showed that the torsion subgroup of  $H^3(X,\Z)$ is nontrivial, and is a birational invariant.
\smallskip	

These three papers have been extremely influential. Though they appeared around the same time, they use very different ideas; in fact, as we will see, the methods tend to apply to different types of varieties. They have been developed and extended, and applied to a number of interesting examples. 
Each of them has its advantages and its drawbacks; very roughly:\smallskip

$\bullet$ The intermediate Jacobian method is quite efficient, but applies only in dimension 3;

$\bullet$ The computation of birational automorphisms leads to the important notion of \emph{birational rigidity}. However it is not easy to work out; so far it applies essentially to Fano varieties of index 1 (see \ref{Fano}), which are not known to be unirational in dimension $>3$.

$\bullet$ Torsion in $H^3$  gives an obstruction to a property weaker than rationality, called  \emph{stable rationality} (\S\ref{stab}). Unfortunately it applies only to very particular varieties, and not to  the standard examples of unirational varieties, like hypersurfaces or complete intersections. However we will discuss in \S \ref{chow} a new idea of C. Voisin which extends considerably the range of that method.

\smallskip	
They are still essentially the basic methods to prove non-rationality results. A notable exception is the method of Koll\'ar using reduction modulo $p$; however it applies only to rather specific examples, which are not known to be unirational. We will describe briefly his results in (\ref{K}).

\smallskip	
A final remark : at the time they were discovered the three methods used the difficult   \emph{resolution of indeterminacies} due to Hironaka. This is a good reason why the Italian algebraic geometers could not succeed! It was later realized  that the birational invariance of $\tor H^3(V,\Z)$ can be proved without appealing to the resolution of singularities, see (\ref{bnr}) -- but this still requires some highly nontrivial algebraic apparatus.

\medskip	
\section{The candidates}
In this section we will introduce various classes of  varieties which are natural candidates to be counter-examples to the L\"uroth problem.

\subsection{Rationality and unirationality}
Let us first recall the basic definitions which appear in the L\"uroth problem. We work over the complex numbers. A \emph{variety} is an integral scheme of finite type over $\C$.
\begin{defi}
$1)$ A variety $V$ is \emph{unirational} if there exists a dominant rational map $\P^n\dasharrow V$.

$2)$ $V$ is \emph{rational} if there exists a birational map $\P^n\bir V$.
\end{defi}

In the definition of unirationality we can take $n=\dim V$: indeed, if we have a dominant rational map $\P^N\dasharrow V$, its restriction to a general linear subspace of dimension $\dim(V)$ is still dominant.

We may rephrase these definitions in terms of the function field $\C(V)$ of $V$ : 
$V$ is unirational if $\C(V)$ is contained in a purely transcendental extension of $\C$; $V$
is rational if $\C(V)$ is a purely transcendental extension of $\C$. Thus the L\"uroth problem asks whether every extension of $\C$ contained in $\C(t_1,\ldots ,t_n )$ is purely transcendental.

\subsection{Rational connectivity}
Though the notion of unirationality is quite natural, it is very difficult to handle. The crucial problem is that so far there is no known method to prove non-unirationality, like the
ones we mentioned in (\ref{mod}) for non-rationality.

There is a  weaker notion which behaves much better than unirationality, and which covers all varieties we will be interested in :
\begin{defi}
A smooth projective variety $V$ is \emph{rationally connected} (\emph{RC} for short) if any two points of $V$ can be joined by a rational curve.
\end{defi}
It is enough to ask that two general points of $V$ can be joined by a rational curve, or even by a chain of rational curves. In particular, rational connectivity is a birational property.

In contrast to unirationality, rational connectivity has extremely good properties (see for instance  \cite{Ar} for proofs and references) :

\smallskip	
$a)$ It is an open and closed property; that is, given a smooth projective morphism $f:V\rightarrow B$ with $B$ connected, if some fiber of $f$ is RC,  \emph{all} the fibers are RC.

$b)$ Let $f:V\dasharrow B$ be a rational dominant map. If $B$ and the general fibers of $f$ are RC, $V$ is RC.

$c)$ If $V$ is RC, all contravariant tensor fields vanish; that is, $H^0(V,(\Omega ^1_V)^{\pt n})=0$ for all $n$. It is conjectured that the converse holds; this is proved in dimension $\leq 3$.

\smallskip	
Neither $a)$ nor $b)$ are expected to hold when we replace rational connectivity by unirationality or rationality.  For $a)$, it is expected that the general quartic threefold is \emph{not} unirational (see \cite[V.9]{Roth}), though some particular $V_4$ are; so unirationality should not be stable under deformation.  Similarly it is expected that the general cubic fourfold is not rational, though some of them are known to be rational. 

Projecting a cubic threefold $V_3$ from a line contained in $V_3$ gives a rational dominant map to $\P^2 $ whose generic fiber is a rational curve, so $b)$ does not hold for rationality.
The same property holds more generally for a general hypersurface of degree $d$ in $\P^4$ with a $(d-2)$-uple  line; it is expected that it is not even unirational for $d\geq 5$ \cite[IV.6]{Roth}.

\subsection{Fano manifolds}\label{Fano}
A more restricted class than RC varieties is that of \emph{Fano manifolds} -- which were extensively studied by Fano  in dimension 3. A smooth projective variety $V$ is \emph{Fano} if the anticanonical bundle $K_V^{-1}$ is ample. This implies that $V$ is RC; but  contrary to the notions considered so far, this is not a property of the birational class of $V$.

A Fano variety $V$ is called \emph{prime} if $\Pic(V)=\Z$ (the classical terminology is ``of the first species"). In that case we have $K_V=L^{-r}$, where $L$ is the positive generator of $\Pic(V)$. The integer $r$ is called the \emph{index} of $V$. Prime Fano varieties are somehow minimal among RC varieties : they do not admit a Mori type 
contraction or morphisms to smaller-dimensional varieties.

In the following table we list what is known about rationality issues for prime Fano threefolds, using their classification by Iskovskikh \cite{IFano} : for each of them, whether it is unirational or rational, and, if it is not rational, the method of proof and the corresponding reference.
 The only Fano threefolds of index $\geq 3$ are $\P^3$ and the smooth quadric $V_2\subset \P^4$, so we start with index 2, then 1:

\bigskip	
\begin{center}\setlength{\tabcolsep}{7pt}

\renewcommand{\arraystretch}{1.5} 
\begin{tabular}{|c|c|c|c|c|}
\hline
variety & unirational & rational & method & reference\\
\hhline{|=====|}

$V_6\subset \P(1,1,1,2,3)$ & ? & no & $\mathrm{Bir}(V)$ & \cite{Gri}\\
\hline
quartic double $\P^3$ & yes & no & $JV$ & \cite{V1}\\
\hline
$V_3\subset \P^4$ & " & no & $JV$ & \cite{CG}\\
\hline 
$V_{2,2}\subset \P^5\ ,\  X_5\subset \P^6$ & " & yes & &\\
\hhline{|=====|}
sextic double $\P^3$ & ? & no & $\mathrm{Bir}(V)$ & \cite{IM}\\
\hline
$ V_4\subset \P^4 $ & some & no & $\mathrm{Bir}(V)$ &\cite{IM} \\
\hline
$ V_{2,3}\subset \P^5 $ & yes & no (generic)&$JV $, $\mathrm{Bir}(V)$ & \cite{Bintjac,Puk}\\ 
\hline
$ V_{2,2,2}\subset \P^6 $ & " &no & $JV$  &\cite{Bintjac}\\
\hline
$ X_{10}\subset \P^7 $ & " &no (generic)& $JV$  &\cite{Bintjac}\\ 
\hline
$ X_{12},X_{16},X_{18},X_{22} $ & " &yes & &\\
\hline
$ X_{14}\subset \P^9 $ & " & no & $JV $& \cite{CG} + \cite{F3}\footnotemark\\
\hline
\end{tabular}
\footnotetext{ Fano  proved in \cite{F3} that the variety $X_{14}$ is birational to a smooth cubic threefold.}
\end{center}

\bigskip	
A few words about notation : as before $V_{d_1,\ldots ,d_p}$ denotes a smooth complete intersection of multidegree $(d_1,\ldots ,d_p)$ in $\P^{p+3}$, or, for the first row, in the weighted projective space $\P(1,1,1,2,3)$. A quartic (resp.\ sextic) double $\P^3$ is a double cover of $\P^3$ branched along a smooth quartic (resp.\ sextic) surface.
The notation $X_d\subset \P^m$  means a smooth threefold of degree $d$ in $\P^m$.
The mention ``(generic)'' means that non-rationality is known only for those varieties belonging to a certain Zariski open subset of the moduli space. 

\subsection{Linear quotients}\label{V/G}
An important source of unirational varieties is provided by the quotients $V/G$, where $G$ is an algebraic group (possibly finite) acting linearly on the vector space $V$. 
These varieties, and the question whether they are rational or not, appear naturally in various situations. The case $G$ finite is known as the Noether problem (over $\C$); we will see below (\ref{bnr}) that a counter-example has been given by Saltman \cite{Sa}, using an elaboration of the Artin-Mumford method. The case where $G$ is a connected linear group  appears in a number of moduli problems, 
but there is still no example where the quotient $V/G$ is known to be non-rational -- in fact the general expectation is that all these quotients should be rational, but this seems out of reach at the moment. 

A typical case is  the moduli space $\mathcal{H}_{d,n}$ of hypersurfaces of degree $d\geq 3$ in $\P^n$, which is birational to $H^0(\P^n,\mathcal{O}_{\P^n}(d))/\mathrm{GL}_{n+1}$ -- more precisely, it is the quotient of the open subset  of forms defining a smooth hypersurface. For $n=2$ the rationality is now known except for a few small values of $d$,  see for instance \cite{BBK} for an up-to-date summary; for $n\geq 3$ there are only a few cases where $\mathcal{H}_{d,n}$ is known to be rational. 
We refer to \cite{Dol} for a survey of results and problems, and to \cite{CTS} for a more recent text.

\medskip	
\section{The intermediate Jacobian}
In this section we discuss our first non-rationality criterion, using the intermediate Jacobian. Then we will give an easy example of a cubic threefold which satisfies this criterion, hence gives a counter-example to the L\"uroth problem.

\subsection{The Clemens-Griffiths criterion}
In order to define the intermediate Jacobian, let us first recall the Hodge-theoretic construction of   the Jacobian of a (smooth, projective) curve $C$. We start from the Hodge decomposition
\[H^1(C,\Z)\subset H^1(C,\C)=H^{1,0}\oplus H^{0,1} \]with $H^{0,1}=\overline{H^{1,0}}$. The latter condition implies that the projection $H^1(C,\R)\rightarrow H^{0,1}$ is a ($\R$-linear) isomorphism, hence that 
the image $\Gamma $ of $H^1(C,\Z)$ in $H^{0,1}$ is a \emph{lattice}  (that is, any basis of $\Gamma $ is a basis of $H^{0,1}$ over $\R$). The quotient $JC:=H^{0,1}/\Gamma $ is a complex torus. But we have more structure. For $\alpha ,\beta \in H^{0,1}$, put
$H(\alpha ,\beta )=2i\int_C \bar{\alpha }\wedge \beta $. Then $H$ is a positive hermitian form on $H^{0,1}$, and the restriction of $\mathrm{Im}(H)$ to $\Gamma \cong H^1(C,\Z)$ coincides with the cup-product
\[H^1(C,\Z)\otimes H^1(C,\Z)\rightarrow H^2(C,\Z)=\Z\ ;\]thus it induces on $\Gamma $ a skew-symmetric, integer-valued form, which is moreover
\emph{unimodular}. In other words, $H$ is a \emph{principal polarization} on $JC$ (see \cite{BL}, or \cite{Ts} for an elementary treatment). This is equivalent to the data of an ample divisor $\Theta \subset JC$ (defined  up to translation) satisfying $\dim H^0(JC,\O_{JC}(\Theta ))=1$. Thus  $(JC,\Theta )$ is a \emph{principally polarized abelian variety} (p.p.a.v. for short), called the \emph{Jacobian} of $C$.

One can mimic this definition for higher dimensional varieties, starting from the odd degree cohomology; this defines the general notion of intermediate Jacobian. In general it is  only a complex torus, not an abelian variety. But the situation is much nicer in the case of interest for us, namely rationally connected threefolds. For such a threefold $V$ we have $H^{3,0}(V)=H^0(V,K_V)=0$, hence the Hodge decomposition for $H^3$ becomes :
\[H^3(V,\Z)_{\mathrm{tf}}\subset H^3(V,\C)=H^{2,1}\oplus H^{1,2} \]with $H^{1,2}=\overline{H^{2,1}}$ ($H^3(V,\Z)_{\mathrm{tf}}$ denotes the quotient of $H^3(V,\Z)$ by its torsion subgroup). As above $H^{1,2}/H^3(V,\Z)_{\mathrm{tf}}$ is a complex torus, with a principal polarization defined by the hermitian form $(\alpha ,\beta )\mapsto -2i\int_V\bar{\alpha }\wedge \beta $: this is the \emph{intermediate Jacobian} $JV$ of $V$.

\medskip	
We will use several times the following well-known and easy lemma, see for instance \cite[ Thm.\ 7.31]{V2} :
\begin{lem}\label{eclat}
Let $X$ be a complex manifold, $Y\subset X$ a closed submanifold of codimension $c$, $\hat{X}$ the variety obtained by blowing up $X$ along $Y$. There are natural isomorphisms  \[ H^p(\hat{X},\Z)\iso  H^p(X,\Z)\oplus \sum\limits_{k=1}^{c-1}H^{p-2k}(Y,\Z)\ .\]
\end{lem}
\begin{thm}[Clemens-Griffiths criterion]\label{CG}
Let $V$ be a smooth rational projective threefold.  The intermediate Jacobian  $JV$ is isomorphic (as p.p.a.v.) to the Jacobian of a curve or to a product of Jacobians.
\end{thm}
\noindent\emph{Sketch of proof} : Let $\varphi :\P^3\bir V$ be a birational map. Hironaka's resolution of indeterminacies provides us with a commutative diagram
\[ \xymatrix{& P\ar[dl]_b\ar[dr]^f &\\  \P^3\ar@{-->}[rr]^\varphi & & V
} \]where $b: P\rightarrow \P^3$ is a composition of blowing up, either of points or of smooth curves, and $f$ is a birational \emph{morphism}. 

We claim that $JP$ is a product of Jacobians of curves. Indeed by  Lemma \ref{eclat},
blowing up a point in a threefold $V$ does not change $H^3(V,\Z)$, hence does not change $JV$ either. If we blow up a smooth curve $C\subset V$ to get a variety $\hat{V}$, Lemma \ref{eclat} gives a canonical isomorphism $H^3(\hat{V},\Z)\cong \allowbreak H^3(V,\Z)\,\oplus \,H^1(C,\Z)$, compatible in an appropriate sense with the Hodge decomposition and the cup-products; this implies $J\hat{V}\cong JV\times JC$ as p.p.a.v. Thus going back to our diagram, we see that $JP$ is isomorphic to $JC_1\times \ldots \times JC_p$, where $C_1,\ldots ,C_p$ are the (smooth) curves which we have blown up in the process.

How do we go back to $JV$? Now we have a birational \emph{morphism} $f:P\rightarrow V$, so we have homomorphisms $f^*: H^3(V,\Z)\rightarrow H^3(P,\Z)$ and $f_*:H^3(P,\Z)\rightarrow H^3(V,\Z)$ with $f_*f^*=1$, again compatible  with the Hodge decomposition and the cup-products in an appropriate sense. Thus $H^3(V,\Z)$, with its polarized Hodge structure, is a direct factor of  $H^3(P,\Z)$; this implies that
 $JV$ is a direct factor of $JP\cong JC_1\times \ldots \times JC_p$, in other words there exists a \pp $A$ such that $JV\times A\cong  JC_1\times \ldots \times JC_p$.
 
 How can we conclude? In most categories the decomposition of an object as a product is not unique (think of vector spaces!). However here a miracle occurs. Let us say that a p.p.a.v. is \emph{indecomposable} if it is not isomorphic to a product of nontrivial p.p.a.v.
 
 \begin{lem}\label{ppav}
$1)$ A p.p.a.v.$\,(A,\Theta )$ is indecomposable if and only if the divisor $\Theta $ is irreducible.

$2)$ Any p.p.a.v. admits a \emph{unique} decomposition as a product of indecomposable p.p.a.v.
\end{lem}
\noindent\textit{Sketch of proof} : We start by recalling some classical properties of abelian varieties, for which we refer to \cite{Mumford}. Let $D$ be a divisor on an abelian variety $A$; for $a\in A$ we denote by $D_a$ the translated divisor $D+a$. The map $\varphi _D:a\mapsto \O_A(D_a-D)$ is a homomorphism from $A$ into its dual variety $\hat{A}$, which parametrizes topologically trivial line bundles on $A$. 
If $D$ defines a principal polarization, this map is an isomorphism. 

Now suppose our \pp $(A,\Theta )$ is a product $(A_1,\Theta _1)\times \ldots \times (A_p,\Theta _p)$. Then $\Theta =\Theta ^{(1)}+\ldots +\Theta ^{(p)}$, with $\Theta ^{(i)}:=A_1\times \ldots \Theta _i\times \ldots \times A_p$;  we recover
 the summand $A_i\subset A$ as $\varphi _{\Theta }^{-1}(\varphi _{\Theta ^{(i)}}(A))$.
 Conversely, let  $(A,\Theta )$ be a p.p.a.v., and let $\Theta ^{(1)},\ldots ,\Theta ^{(p)}$ be the irreducible components of $\Theta $ (each of them occurs with multiplicity one, since otherwise one would have $h^0(A;\O_A(\Theta ))>1$). Putting $A_i:=\varphi _{\Theta }^{-1}(\varphi _{\Theta ^{(i)}}(A))$ and $\Theta _i:=\Theta ^{(i)}_{|A_i}$, it is not difficult to check that $(A,\Theta )$ is the product of the $(A_i,\Theta _i)$ -- see \cite{CG}, Lemma 3.20 for the details.\qed
 
 \medskip	
Once we have this, we conclude as follows. The Theta divisor of a Jacobian $JC$ is the image of the Abel-Jacobi map $C^{(g-1)}\rightarrow JC$, and therefore is irreducible. From the isomorphism $JV\times A\cong  \allowbreak JC_1\times \ldots \times JC_p$ and the Lemma we conclude that $JV$ is isomorphic to $JC_{i_1}\times \ldots \times JC_{i_r}$ for some subset $\{i_1,\ldots ,i_r\} $ of $[1,p]$.\qed

\medskip	
\noindent\emph{Remark}$.-$
One might think that  products of Jacobians are more general than Jacobians, but it goes the other way around: in the moduli space $\mathcal{A}_g$ of $g$-dimensional p.p.a.v., the boundary $\bar{\mathcal{J}_g}\smallsetminus \mathcal{J}_g$ of the Jacobian locus 
 is precisely the locus of products of lower-dimensional Jacobians.

\subsection{The Schottky problem}\label{schott}

Thus to show that a threefold $V$ is not rational, it suffices to prove that its intermediate Jacobian is not the Jacobian of a curve, or a product of Jacobians. Here we come across  the  classical Schottky problem : the characterization of Jacobians among all p.p.a.v. (the usual formulation of the Schottky problem asks  for equations of the Jacobian locus inside the moduli space of p.p.a.v.; here we are more interested in special geometric properties of Jacobians). One frequently used approach is through the singularities of the Theta divisor : the dimension of $\Sing(\Theta )$ is $\geq g-4$ 
for a Jacobian $(JC,\Theta )$ of dimension $g$, and  $g-2$ for a product. However controlling $\Sing(\Theta )$ for an intermediate Jacobian is quite difficult, and requires a lot of information on the geometry of $V$. Let us just give a sample :
\begin{thm}
Let $V_3\subset \P^4$ be a smooth cubic threefold. The divisor $\Theta \subset JV_3$ has a unique singular point $p$, which is a triple point. The tangent cone $\P T_p(\Theta )\subset \P T_p(JV_3)\cong \P^4$ is isomorphic to $V_3$.
\end{thm}
This elegant result, apparently due to Mumford (see \cite{Bcubic} for a proof), implies both the non-rationality of $V_3$ (because $\dim \Sing(\Theta )=0$ and $\dim JV_3=5$) and the \emph{Torelli theorem} : the cubic $V_3$ can be recovered from its (polarized) intermediate Jacobian. 

There are actually few cases where we can control so well the singular locus of  the Theta divisor. One of these is the quartic double solid, for which $\Sing(\Theta )$ has a component of codimension 5 in $JV$ \cite{V1}. Another case is that of  \emph{conic bundles}, that is, threefolds $V$ with a flat morphism $p:V\rightarrow \P^2$, such that for each closed point $s\in \P^2$ the fiber $p^{-1}(s)$ is isomorphic to a plane conic (possibly singular). In that case $JV$ is a \emph{Prym variety}, associated to a natural double covering of the \emph{discriminant curve} $\Delta \subset\P^2$ (the locus of $s\in \P^2$ such that $p^{-1}(s)$ is singular). Thanks to Mumford we have some control on the singularities of  the Theta divisor of a Prym variety, enough to show that $JV$ \emph{is not a Jacobian (or a product of Jacobians) if $\deg(\Delta )\geq 6$} \cite[thm.\ 4.9]{Bintjac}. 

Unfortunately, apart from the cubic,  the only prime Fano threefold to which this result applies is the $V_{2,2,2}$ in $\P^6$. However, 
the Clemens-Griffiths criterion of non-rationality is an \emph{open} condition. In fact, we have a stronger result,  which follows  from the properties of the \emph{Satake compactification}  of $\mathcal{A}_g$ \cite[lemme 5.6.1]{Bintjac} :
 \begin{lem}\label{defjac}
Let $\pi :V\rightarrow B$ be a flat  family of projective threefolds over a smooth curve $B$. Let $\mathrm{o}\in B$; assume that :

$\bullet$ The fiber $V_b$ is smooth for $b\neq \mathrm{o}$;

 $\bullet$ $V_{\mathrm{o}}$ has only ordinary double points; 
 
 $\bullet$ For a desingularization $\tilde{V}_{\mathrm{o}} $ of $V_{\mathrm{o}}$, 
 $J\tilde{V}_{\mathrm{o}}$ is not a Jacobian or a product of Jacobians.
 
 Then for $b$ outside a finite subset   of $B$,   $V_b$ is not rational.
\end{lem}
From this we deduce the generic non-rationality statements of  (\ref{Fano}) \cite[Thm.\ 5.6]{Bintjac} : in each case one finds a degeneration as in the Lemma, such that $\tilde{V}_{\mathrm{o}} $ is a conic bundle with a discriminant curve of degree $\geq 6$, hence the Lemma applies.

\subsection{An easy counter-example} 
The results of the previous section require rather involved methods. We will now discuss a much more elementary approach, which unfortunately applies only to specific varieties.

\begin{thm}
The  cubic threefold $V\subset \P^4$ defined by $\sum\limits_{i\in \Z/5}X_i^2X_{i+1}=0$ is not rational.
\end{thm}

\pr Let us first prove that $JV$ is not a Jacobian. Let $\zeta  $ be a primitive $11$-th root of unity.
The key point is that $V$ admits the automorphisms\par
$\delta : (X_0,X_1,X_2,X_3,X_4 )\mapsto (X_0,\zeta X_1,\zeta ^{-1}X_2,\zeta ^3X_3,\zeta ^{6}X_4)$,\par
$\sigma : (X_0,X_1,X_2,X_3,X_4 )\mapsto (X_1,X_2,X_3,X_4,X_0 )$,\par
\noindent which satisfy
 $\delta ^{11}=\sigma ^5=1$ and $\sigma \delta \sigma ^{-1}= \delta ^{-2}$. 

They induce automorphisms $\delta ^*,\sigma ^*$ of $JV$. 
Suppose that $JV$ is isomorphic (as p.p.a.v.) to the Jacobian $JC$ of a curve $C$. The Torelli theorem for curves gives an exact sequence
\[1\rightarrow \Aut(C)\rightarrow \Aut(JC)\rightarrow \Z/2\ ;\]since  $\delta ^*$ and  
 $\sigma ^*$ have odd order, they are induced by
 automorphisms $\delta_C,\sigma_C$ of $C$, satisfying  $\sigma_C\delta_C\sigma_C^{-1}=\delta_C^{-2}$.

Now we apply the Lefschetz fixed point formula. The automorphism $\delta $ of $V$ fixes the 5 points corresponding to the basis vectors of $\C^5$; it acts trivially on $H^{2i}(V,\Q)$ for $i=0,\ldots ,3$. Therefore we find $\Tr \delta ^*_{\ |H^3(V,\Q)}=-5+4=-1$. Similarly $\sigma $ fixes the 4 points $(1,\alpha ,\alpha ^2,\alpha ^3,\alpha ^4)$ of $V$ with $\alpha ^5=1$, $\alpha \neq 1$, so $\Tr \sigma ^*_{\ |H^3(V,\Q)}=-4+4=0$. 

Applying now the Lefschetz  formula to $C$, we find that $\sigma_C$ has two fixed points on $C$ and $\delta_C$ three. But since $\sigma_C$ normalizes the subgroup generated by $\delta_C$, it preserves the 3-points set $\mathrm{Fix}(\delta_C)$; since it is of order 5, it must fix each of these 3 points, which gives a contradiction.

Finally suppose $JV$ is isomorphic to a product $A_1\times \ldots \times A_p$ of p.p.a.v. By the unicity lemma (Lemma \ref{ppav}), the automorphism $\delta ^*$ permutes the factors $A_i$. Since $\delta $ has order 11 and $p\leq 5$, this permutation must be trivial, so $\delta ^*$ induces an automorphism of  $A_i$ for each $i$, hence of $H^1(A_i,\Q)$; but the group $\Z/11$ has only one nontrivial irreducible representation defined over $\Q$,  given by the cyclotomic field $\Q(\zeta )$, with $[\Q(\zeta ):\Q]=10$. Since $\dim (A_i)<5$ we see that the action of $\delta ^*$ on each $A_i$, and therefore on $JV$, is trivial. But this contradicts the relation $\Tr \delta ^*_{\ |H^3(V,\Q)}=-1$.\qed

\medskip	
\noindent\emph{Remarks}$.-$ 1) The cubic $V$ is the \emph{Klein cubic threefold}; it is birational to the moduli space of abelian surfaces with a polarization of type $(1,11)$ \cite{GP}. In particular it admits an action of the group $\mathrm{PSL}_2(\F_{11})$ of order 660, which is in fact its automorphism group
 \cite{Adler}. From this one could immediately conclude by using the Hurwitz bound $\#\Aut(C)\leq 84(g(C)-1)$ (see \cite{Bsex}). 
 
2) This method applies to other threefolds for which the non-rationality was not previously known, in particular the $\mathfrak{S}_7$-symmetric $V_{2,3}$ given by $\sum X_i=\sum X_i^2=\sum X_i^3=0$ in $\P^6$ \cite{Bsex} or the $\mathfrak{S}_6$-symmetric $V_{4}$ with 30 nodes given by $ \sum X_i=\sum X_i^4=0 $ in $\P^5 $ \cite{Bquar}.

\medskip
\section{Two other methods}
In this section we will briefly present  two other ways to get non-rationality results for certain Fano varieties. Let us stress that in dimension $\geq 4$ these varieties are not known to be unirational, so these methods do not give us new counter-examples to the L\"uroth problem.

\subsection{Birational rigidity}
As  mentioned in the introduction, Iskovskikh and Manin proved that a smooth quartic threefold $V_4\subset \P^4$ is not rational by proving that any birational automorphism of $V_4$ is actually biregular. But they proved much more, namely that $V_4$ is \emph{birationally superrigid} in the following sense :

\begin{defi}
Let $V$ be a prime Fano variety (\ref{Fano}). We say that $V$ is \emph{birationally rigid} if :

$a)$ There is no rational dominant map $V\dasharrow S$ with $0<\dim(S)<\dim(V)$ and with general fibers of  Kodaira dimension $-\infty$;

$b) $ If $V$ is birational to another prime Fano variety $W$, then $V$ is isomorphic to $W$. 

\noindent We say that $V$ is \emph{birationally superrigid} if any birational map $V\bir W$ as in $b)$ is an isomorphism.
\end{defi}
(The variety $W$ in $b)$ is allowed to have certain mild singularities, the so-called $\Q$-factorial terminal singularities.) 

After the  pioneering work \cite{IM}, birational (super)rigidity has been proved for a number of Fano varieties of index 1. 
Here is a sample; we refer to the surveys \cite{Puk} and \cite{Chel} for ideas  of proofs and for many more examples.

\smallskip	
$\bullet$ Any smooth hypersurface of degree $n$ in $\P^n$ is  birationally superrigid \cite{dF}.

$\bullet$ A general $V_{2,3}$ in $\P^5$ is birationally rigid. It is not birationally superrigid, since it contains a curve of lines, and each line defines by projection a 2-to-1 map to $\P^3$, hence a birational involution of $V_{2,3}$.

$\bullet$ A general $V_{d_1,\ldots ,d_c}$ in $\P^n$ of index 1 (that is, $\sum d_i=n$) with $n> 3c$  is birationally superrigid. 

$\bullet$ A double cover of $\P^n$ branched along a smooth hypersurface of degree $2n$ is birationally superrigid.

\subsection{Reduction to characteristic $p$}\label{K}
\begin{thm}\cite{K}
For $d\geq 2\lceil \dfrac{n+3}{3}\rceil  $, a very general hypersurface $V_d\subset \P^{n+1}$ is not ruled, and in particular not rational. 
\end{thm}
A variety is \emph{ruled} if it is birational to $W\times \P^1$ for some variety $W$. 
  ``Very general" means that the corresponding point in the space parametrizing our hypersurfaces lies outside a countable union of strict closed subvarieties. 
  
  \medskip	
 The bound $d\geq 2\lceil \dfrac{n+3}{3}\rceil  $ has been  lowered to  $d\geq 2\lceil \dfrac{n+2}{3}\rceil  $ by Totaro \cite{T}; this implies in particular that a very general quartic fourfold is not rational. More important, by combining Koll\'ar's method with a new idea of Voisin (see \S\ref{chow}), Totaro shows that a  very general  $V_d\subset \P^{n+1}$ with $d$ as above is  not \emph{stably rational} (\S\ref{stab}).

\smallskip	
Let us give a very rough idea of Koll\'ar's proof, in the case $d$ is even. It starts from the well-known fact that the hypersurface  $V_d$   specializes to  a double covering $Y$ of a hypersurface of degree $d/2$. This can be still done in characteristic 2, at the price of getting some singularities on $Y$, which must be resolved. The reward is that the resolution $Y'$ of $Y$ has a very unstable tangent bundle; more precisely, $\Omega ^{n-1}_{Y'}$ contains a positive line bundle, and this prevents $Y'$ to be ruled. Then a general result of Matsusaka implies that a very general $V_d$ cannot be ruled.

\medskip	
\section{Stable rationality}\label{stab}
There is an intermediate notion between rationality and unirationality which turns out to be important :
\begin{defi}
A variety $V$ is \emph{stably rational} if $V\times \P^n$ is rational for some $n\geq 0$.
\end{defi}
In terms of field theory, this means that $\C(V)(t_1,\ldots ,t_n)$ is a purely transcendantal extension of $\C$.

Clearly, rational $\ \Rightarrow\ $ stably rational $\ \Rightarrow\ $ unirational. We will see that these implications are strict. For the first one, we have :
\begin{thm}\cite{BCTSSD}
Let $P(x,t)=x^3+p(t)x+q(t)$  be an irreducible polynomial in $\C[x,t]$, whose discriminant 
$\delta (t):=4p(t)^3+27q(t)^2$ has degree $\geq 5$. The affine hypersurface  $V\subset \C^4$ defined by 
$y^2-\delta (t)z^2=P(x,t)$ is stably rational but not rational.
\end{thm}
This answered a  question asked by Zariski in 1949 \cite{SgZar}.

The non-rationality of $V$ is proved using the intermediate Jacobian, which turns out to be the Prym variety associated to an admissible double covering of nodal curves. The stable rationality, more precisely the fact that $V\times \P^3$ is rational, was proved in \cite{BCTSSD} using some particular torsors under certain algebraic tori. A different construction due to Shepherd-Barron shows that actually $V\times \P^2$ is rational \cite{ShB}; we do not know whether 
$V\times \P^1$ is rational.

\medskip	
To find unirational varieties which are not stably rational, we cannot use the Clemens-Griffiths criterion since it applies only in dimension 3. The group of birational  automorphisms is very complicated  for a variety of the form $V\times \P^n$; so the only available method is the torsion of $H^3(V,\Z)$ and its subsequent refinements, which we will examine in the next sections.

\medskip	
\noindent\emph{Remark}$.-$\label{retract} There are other notions lying between unirationality and rationality. Let us say that a variety $V$ is

$\bullet$ retract rational if there exists a rational dominant map $\P^N\dasharrow V$ which admits a rational section;

$\bullet$ factor-rational if there exists another variety $V'$ such that $V\times V'$ is rational.

We have  the implications :

\centerline{rational $\ \Rightarrow\ $ stably rational $\ \Rightarrow\ $ factor-rational $\ \Rightarrow\ $ retract rational $\ \Rightarrow\ $ unirational.}

Unfortunately at the moment we have no examples (even conjectural) of varieties which are retract rational but not stably rational. For this reason we will focus on the stable rationality, which seems at this time the most useful of these notions. Indeed we will see now that there are some classes of linear quotients  $V/G$ (see \ref{V/G}) for which  we can prove stable rationality.

\medskip	
Let $G$ be a  reductive group acting on a variety $V$. We   say that the action is \emph{almost free} if there is a nonempty Zariski open subset $U$ of $V$ such that the stabilizer of each point of $U$ is trivial. 

\begin{prop}
Suppose that there exists an almost free linear representation $V$ of $G$ such that the quotient $V/G$ is rational. Then for every almost free representation $W$ of $G$, the quotient $W/G$ is stably rational. 
\end{prop}
The proof goes as follows \cite{Dol} : let $V^{\mathrm{o}}$ be a  Zariski open subset  of $V$ where $G$ acts freely. Consider the diagonal action of $G$ on $V^{\mathrm{o}}\times W$;  standard arguments (the ``no-name lemma") show that  the projection $(V^{\mathrm{o}}\times W)/G\rightarrow V^{\mathrm{o}}/G$ defines a vector bundle over $V^{\mathrm{o}}/G$. Thus $(V\times W)/G$ is birational to $(V/G)\times W$ (which is a rational variety), and  symmetrically to $V\times (W/G)$, so $W/G$ is stably rational.\qed

\smallskip	
For many groups it is easy to find an almost free representation with rational quotient : this is the case for instance for a subgroup $G$ of $\mathrm{GL}_n$ such that the quotient $\mathrm{GL_n}/G$ is rational (use the linear action of $\mathrm{GL}_n$ on $\mathrm{M}_n(\C)$ by multiplication). This applies to 
$\mathrm{GL}_n$, $\mathrm{SL}_n$, $\mathrm{O}_n$ ($\mathrm{GL}_n/\mathrm{O}_n$ is  the space of non-degenerate quadratic forms), $\mathrm{SO}_n$, $\mathrm{Sp}_n$ etc. In \cite{Bogo} Bogomolov proves that every
  simply connected simple group has this property, except perhaps  $E_8$. 
  
\medskip	 
This gives many examples of stably rational varieties. For instance, 
 the moduli space $\mathcal{H}_{d,n}$ of hypersurfaces of degree $d$ in $\P^n$ (\ref{V/G})  is stably rational when $d\equiv 1\ \mathrm{mod.}\ (n+1)$  : the standard representation $\rho $ of $\mathrm{GL}_{n+1}$ on $H^0(\P^n,\O_{\P^n}(d))$ is not almost free, 
 but  the representation $\rho \otimes \det^k$, with $k=\frac{1-d}{n+1} $, is almost free and gives the same quotient.
 
 \medskip	
\section{The torsion of $H^3(V,\Z)$ and the Brauer group}
\subsection{Birational invariance}
Artin and Mumford used the following property of stably rational varieties :
\begin{prop}\label{am}
Let $V$ be a stably rational projective manifold. Then $H^3(V,\Z)$ is torsion free.
\end{prop}
\pr The K\"unneth formula gives an isomorphism  $H^3(V\times \P^m,\Z)\cong H^3(V,\Z)\,\oplus\, H^1(V,\Z)$; since $H^1(V,\Z)$ is torsion free the torsion subgroups of $H^3(V,\Z)$ and $H^3(V\times \P^m,\Z)$ are isomorphic,
 hence  replacing $V$ by $V\times \P^m$ we may assume that $V$ is rational.
Let $\varphi :\P^n\bir V$ be a birational map. As in the proof of the Clemens-Griffiths criterion, we have Hironaka's ``little roof''
\[ \xymatrix{& P\ar[dl]_b\ar[dr]^f &\\  \P^n\ar@{-->}[rr]^\varphi & & V
} \]where $b: P\rightarrow \P^n$ is a composition of blowing up  of smooth subvarieties, and $f$ is a birational \emph{morphism}. 

By  Lemma \ref{eclat}, we have $H^3(P,\Z)\cong H^1(Y_1,\Z)\oplus \ldots \oplus H^1(Y_p,\Z)$, where $Y_1,\ldots ,Y_p$ are the subvarieties successively blown up by $b$; therefore   $H^3(P,\Z)$ is torsion free. As in the proof of Theorem \ref{CG},  $H^3(V,\Z)$ is a direct summand of $H^3(P,\Z)$, hence is also torsion free.\qed

\medskip	
We will indicate below (\ref{bnr}) another proof which does not use Hironaka's difficult theorem.

\medskip
\subsection{The Brauer group}\label{brauer}
The torsion of $H^3(V,\Z)$ is strongly related to the \emph{Brauer group} of $V$. There is a huge literature on the Brauer group in algebraic geometry, starting with the three ``expos\'es" by Grothendieck in \cite{G10}. We recall here the cohomological definition(s) of this group; we refer to \cite{G10} for the relation with Azumaya algebras.
\newcommand\br{\operatorname{Br}}
\begin{prop}
Let $V$ be a  smooth  variety. The following definitions are equivalent, and define the \emph{Brauer group} of $V$ :

{\rm (i)} $\br(V)=\Coker c_1:\Pic(V)\otimes \Q/\Z \rightarrow  H^2(V,\Q/\Z) $;

{\rm (ii)} $\br(V)=H^2_{\text{\rm \'et}}(V,\G_m)$ (\'etale cohomology).
\end{prop}
\pr Let $n\in \mathbb{N}$. The exact sequence of \'etale sheaves
$\ 1\rightarrow \Z/n \rightarrow \G_m \qfl{\times n}\G_m\rightarrow 1\ $
gives a cohomology exact sequence 
\[0\rightarrow \Pic(V)\otimes \Z/n \qfl{c_1} H^2(V,\Z/n)\longrightarrow \br(V)\qfl{\times n} \br(V)\ .\](Note that the \'etale cohomology $H^*_{\textrm{\'et}}(V,\Z/n)$ is canonically isomorphic to the classical cohomology). 

Taking the direct limit with respect to  $n$ gives an exact sequence
\begin{equation}\label{pic}
0\rightarrow \Pic(V)\otimes \Q/\Z \qfl{c_1} H^2(V,\Q/\Z)\longrightarrow \tor \br(V)\rightarrow 0\ ;\end{equation}
 it is not difficult to prove that $\br(V)$ is a torsion group  \cite[II, Prop.\ 1.4]{G10}, hence 
the equivalence of the definitions (i) and (ii).\qed

\medskip	
\noindent\emph{Remark}$.-$
If $V$ is compact, the same argument shows that $\br(V)$ is also isomorphic to the torsion subgroup of  $H^2(V,\O_h^*)$, where $\O_h$  is the sheaf of holomorphic functions on $V$ (for the classical topology).

\begin{prop}\label{brH3}
There is a surjective homomorphism $\br(V)\rightarrow \tor H^3(V,\Z)$, which is bijective if $c_1:\Pic(V)\rightarrow H^2(V,\Z)$ is surjective.
\end{prop}
The latter condition is satisfied  in particular if $V$ is projective and $H^2(V,\O_V)=0$.

\medskip	
\pr The exact sequence $0\rightarrow \Z\rightarrow \Q\rightarrow \Q/\Z\rightarrow 0$ gives a cohomology exact sequence
\[0\rightarrow H^2(V,\Z)\otimes \Q/\Z \longrightarrow H^2(V,\Q/\Z)\longrightarrow \tor H^3(V,\Z)\rightarrow 0\ .\]Together  with (\ref{pic}) we get a commutative diagram
\[\xymatrix@M=5pt{0\longrightarrow \Pic(V)\otimes \Q/\Z\ar[r]\ar@<15pt>[d]_{c_1}  & H^2(V,\Q/\Z)\ar[r]\ar@^{=}[d]  & \mathrm{Br}(V)  \longrightarrow  0 \\
0\longrightarrow  H^2(V,\Z)\otimes \Q/\Z\ar[r] & H^2(V,\Q/\Z)\ar[r]& \tor H^3(V,\Z)\longrightarrow 0}\]
 which implies the Proposition.\qed

\medskip	 
We will now describe a geometric way to construct 
 nontrivial elements of the Brauer group.
 \begin{defi}
Let $V$ be a complex variety. A \emph{$\P^m$-bundle} over $V$ is a smooth map $p:P\rightarrow V$ whose geometric fibers are isomorphic to $\P^m$. 
\end{defi}
An obvious example is the projective bundle $\P_V(E)$ associated to a vector bundle $E$ of rank $m+1$ on $V$;   we will actually be interested in those $\P^m$-bundles which are \emph{not} projective.

It is not difficult to see that a $\P^m$-bundle is locally trivial for the \'etale topology. This implies that isomorphism classes of $\P^{n-1}$-bundles over $V$ are parametrized by the \'etale cohomology set $H^1(V, PGL_n)$, where for an algebraic group $\mathrm{G}$ we denote by $G$ the sheaf of local maps to $G$.
The exact sequence of  sheaves of groups
\[1\rightarrow \G_m\rightarrow GL_{n}\rightarrow PGL_{n} \rightarrow 1\]
gives rise to a sequence of pointed sets
\[ H^1(V, GL_{n})\qfl{p} H^1(V, PGL_{n})\qfl{\partial } H^2(V,\G_m)\]
which is exact in the sense that $\partial ^{-1}(1)=\im p$. Thus we associate to each $\P^1$-bundle $p:P\rightarrow V$ a class $[p]$ in $H^2(V,\G_m)$, and this class is trivial if and only if $p$ is  a projective bundle. Moreover, by comparing with the exact sequence $0\rightarrow \Z/n\rightarrow SL_{n}\rightarrow PGL_{n}\rightarrow 1$ we get a commutative diagram
\[\xymatrix{H^1(V,SL_{n})\ar[r]\ar[d] & H^1(V,PGL_{n})\ar[r]\ar@{=}[d]& H^2(V,\Z/n)\ar[d]\\
H^1(V,GL_{n})\ar[r] & H^1(V,PGL_{n})\ar[r]^>>>>{\partial } & H^2(V,\G_m)
}\]
which shows that the image of $\partial $ is contained in the $n$-torsion subgroup of $\br(V)$. 

\medskip	
\subsection{The Artin-Mumford example}

The Artin-Mumford counter-example is a double cover of $\P^3$ branched along a \emph{quartic symmetroid}, that is, a quartic surface defined by the vanishing of a symmetric determinant. 

We start with a web $\Pi $ of quadrics in $\P^3$; its elements  are defined by quadratic forms $\lambda _0q_0+\ldots +\lambda _3q_3$. We  assume that :

(i) $\Pi $ is base point free;

(ii) If a line  in $\P^3$ is singular for a quadric of $\Pi $, it is not contained in another quadric of $\Pi $.

Let $\Delta \subset \Pi $ be the discriminant locus, corresponding to quadrics of rank $\leq 3$. It is a quartic surface (defined by $\det (\sum \lambda _iq_i)=0$);
under our hypotheses, its has 10 ordinary double points, corresponding to quadrics of rank 2, and no other singularity (see for instance \cite{Co}). Let $\pi :V'\rightarrow \Pi $ be the double covering branched along $\Delta $.
Again $V'$ has 10 ordinary double points; blowing up these points we obtain the
 Artin-Mumford threefold $V$. 
 
 Observe that a quadric $q\in\Pi $ has two systems of generatrices (= lines contained in $q$) if $q\in \Pi \smallsetminus \Delta $, and one if $q\in \Delta \smallsetminus \Sing(\Delta )$. Thus the smooth part $V^{\mathrm{o}}$ of $V$ parametrizes pairs $(q, \lambda )$, where $q\in\Pi $ and $\lambda $ is a family of generatrices of $q$.

 \begin{thm}\label{torAM}
The threefold $V$ is unirational but not stably rational.
\end{thm}
\pr 
 Let $\G$ be the Grassmannian of lines in $\P^3$. A general line  is contained in a unique quadric of $\Pi $, and in a unique system of generatrices of this quadric; this defines a dominant rational map $\gamma :\G\dasharrow V' $, thus $V$ is unirational.
We will deduce from  Proposition \ref{am} that $V$ is not stably rational, by proving that $H^3(V,\Z)$ contains an element of order 2. This is done  by a direct calculation in \cite{AM} and, with a different method, in \cite{Bop}; here we will use a more elaborate approach based on the Brauer group.

Consider the variety $P\subset \G\times \Pi $ consisting of pairs $(\ell,q)$ with $\ell\subset q$.  The projection $P\rightarrow \Pi$ factors through a morphism  $p':P\rightarrow V'$. Put $V^{\mathrm{o}}:=V'\smallsetminus \Sing(V')$, and $P^{\mathrm{o}}:=p'^{-1}(V^{\mathrm{o}})$. The restriction $p:P^{\mathrm{o}}\rightarrow V^{\mathrm{o}}$ is a $\P^1$-bundle: a point of $V^{\mathrm{o}}$ is a pair $(q,\sigma )$, where $q$ is a quadric  in $\Pi $ and $\sigma $
a system of generatrices of $q$; the fiber $p^{-1}(q,\sigma )$ is the smooth rational curve parametrizing the lines of $\sigma $. 

\begin{prop}\label{bram}
The $\P^1$-bundle $p:P^{\mathrm{o}}\rightarrow V^{\mathrm{o}}$ does not admit a rational section.\end{prop}
\pr Suppose it does. For a general point $q$ of $\Pi $, 
 the section  maps the two points of $\pi ^{-1}(q)$ to 
two generatrices  of the quadric $q$, one in each system.  These two generatrices intersect in one point $s(q) $ of $q$; the map $q\mapsto s(q)$ is a  rational section of the universal family of quadrics $\mathcal{Q}\rightarrow \Pi $, defined by   $\mathcal{Q}:=\{(q,x)\in  \Pi\times \P^3 \ |\ x\in q \} $.
 This contradicts the following  lemma:
\begin{lem}
Let $\Pi\subset \P(H^0(\P^n, \mathcal{O}_{\P^n}(d))$ be a base point free linear system of hypersurfaces, of degree $d\geq 2$. Consider the universal family $p:\mathcal{H}\rightarrow \Pi $, with $\mathcal{H}:=\{(h,x)\in \Pi \times \P^n\ |\  x\in h\}$. Then $p$ has no rational section.
\end{lem}
\pr Since $\Pi$ is base point free, the second projection $q:\mathcal{H}\rightarrow \P^n$ is a projective bundle, hence $\mathcal{H}$ is smooth. If $p$ has a rational section, the closure $Z\subset \mathcal{H}$ of its image gives a cohomology class $[Z]\in H^{2n-2}(\mathcal{H},\Z)$ such that $p_*([Z])=1$ in $H^0(\Pi ,\Z)$. Let us show that this is impossible.

We have $\dim(\Pi)\geq n$, hence
$2n-2<n-1+\dim(\Pi) =\dim(\mathcal{H})$. By the Lefschetz hyperplane theorem, the restriction map $H^{2n-2}(\Pi \times \P^n,\Z)\rightarrow H^{2n-2}(\mathcal{H},\Z)$ is an isomorphism. Thus $H^{2n-2}(\mathcal{H},\Z)$ is spanned by the classes $p^*h_{\Pi}^i\cdot q^*h_{\P}^{n-1-i}$ for $0\leq i\leq n-1$, where $h_{\Pi}$ and $h_{\P }$ are the hyperplane classes. All these classes go to 0 under $p_*$ except $q^*h_{\P}^{n-1}$, whose degree on each fiber is $d$. Thus the image of $p_*: H^{2n-2}(\mathcal{H},\Z)\rightarrow H^0(\Pi,\Z)=\Z$ is $d\Z$. This proves the lemma, hence the Proposition.\qed 

\medskip	
Thus the $\P^1$-bundle $p$ over $V^{\mathrm{o}}$ is not a projective bundle, hence gives a nonzero 2-torsion class in $\br(V^{\mathrm{o}})$. In the commutative diagram
\[\xymatrix{ \Pic(V )\ar[r]^{c_1}\ar[d] & H^2(V,\Z )\ar[d]^r\ \\
 \Pic(V^{\mathrm{o}}) \ar[r]^{c_1} & H^2(V^{\mathrm{o}},\Z )}\] the top horizontal arrow is surjective because $H^2(V ,\mathcal{O}_{V})=0$. Since $Q:=V\smallsetminus V^{\mathrm{o}}$ is a disjoint union of quadrics, the Gysin exact sequence
  $\ H^2(V ,\Z)\xrightarrow{\,r\,}  H^2(V^{\mathrm{o}},\Z)\rightarrow H^1(Q,\Z)=0 $ shows that $r$ is surjective. Therefore  the map $c_1:\Pic(V^{\mathrm{o}})\rightarrow H^2(V^{\mathrm{o}},\Z)$ is surjective, and by Proposition \ref{brH3} we get a nonzero $2$-torsion class in $H^3(V^{\mathrm{o}},\Z)$. Using again the Gysin exact sequence
\[ 0\rightarrow H^3(V ,\Z)\rightarrow  H^3(V^{\mathrm{o}},\Z)\rightarrow H^2(Q,\Z) \]we find that $\tor H^3(V,\Z)$ is isomorphic to $\tor H^3(V^{\mathrm{o}},\Z)$, hence nonzero.\qed

\medskip	
\subsection{Higher dimension}
The construction of the Artin-Mumford example extends in higher dimension. Let us consider a sufficiently general linear system $\Pi $ of quadrics in $\P^3$, of dimension $n$. We have \[\Sigma \subset \Delta \subset \Pi \]where the discriminant locus $\Delta $ is a quartic hypersurface, and $\Sigma =\mathrm{Sing}(\Delta )$ parametrizes the quadrics of rank $\leq 2$ in $\Pi $; we have $\dim(\Sigma )=n-3$. \emph{We will assume $n\leq 5$}; this guarantees that $\Sigma $ is smooth.

We consider again  the double covering $V'\rightarrow \Pi $ branched along $\Delta $. It is not difficult to see that locally (for the complex topology) around $\Sigma $, $\Delta $ is isomorphic to $\Sigma \times \mathcal{C}$, where $\mathcal{C}$ is a quadratic cone in $\P^3$. It follows that blowing up $V'$ along $\Sigma $ provides a resolution $V\rightarrow V'$.

\begin{prop}\cite{B2Q}\label{2P}
$V$ is unirational, and satisfies $\tor H^3(V,\Z)\neq 0$.
\end{prop}
\emph{Sketch of proof} : Consider again the variety $P\subset \G\times \Pi $ of 
  pairs $(\ell,q)$ with $\ell\subset q$. The projection $P\rightarrow \G$ is a projective bundle over some open subset of $\G$, hence $P$ is rational, and $V$ is unirational.
  
   The construction of the $\P^1$-bundle $P^{\mathrm{o}}\rightarrow V^{\mathrm{o}}$, and the proof that it has no rational section, extend identically. The proof that its class in $H^3(V^{\mathrm{o}},\Z)$ comes from a nonzero class in $H^3(V,\Z)$ requires some more work, for which we refer to \cite{B2Q}.\qed

\medskip	
\subsection{The unramified Brauer group}\label{bnr}
An advantage of the group $\br(V)$ is that 
it can be identified with the
 \emph{unramified Brauer group} $\br_{\textrm{nr}}(\C(V))$, which is defined purely in terms of the field $\C(V)$; this gives directly its birational invariance, without using Hironaka's theorem. Let us explain briefly how this works.

\begin{prop}\label{sexB}
Let $V$ be a smooth projective variety, and $\mathcal{D}$ be the set of integral divisors on $V$. 
For $D$ in $\mathcal{D}$, put $D_{sm}:=D\smallsetminus \Sing(D)$. There is an exact sequence\[0\rightarrow \br(V)\rightarrow \varinjlim_U \br(U)\rightarrow \bigoplus_{D\in\mathcal{D}}H^1(D_{sm},\Q/\Z)\]
where the direct limit is taken over the set of Zariski open subsets  $U\subset V$.
\end{prop}
\pr Let $D$ be an effective reduced divisor on $V$,  and let $U=V\smallsetminus D$. Since $\Sing(D)$ has codimension $\geq 2$ in $V$, the restriction map $H^2(V,\Q/\Z)\rightarrow H^2(V\smallsetminus \Sing(D),\Q/\Z)$ is an isomorphism. Thus we can write part of the Gysin exact sequence as
\[H^0(D_{sm},\Q/\Z)\rightarrow H^2(V,\Q/\Z)\rightarrow H^2(U,\Q/\Z)\rightarrow H^1(D_{sm},\Q/\Z)\ .\]
Comparing with the analogous exact sequence for Picard groups gives a commutative diagram
\[\xymatrix@M=5pt{H^0(D_{sm},\Z)\otimes \Q/\Z \ar[r]\ar[d]^{\wr} &\Pic(V)\otimes \Q/\Z \ar[r]\ar[d]&\Pic(U)\otimes \Q/\Z \ar[r]\ar[d]& 0 &\\
H^0(D_{sm},\Q/\Z) \ar[r] &H^2(V,\Q/\Z) \ar[r]&H^2(U, \Q/\Z )\ar[r]& H^1(D_{sm},\Q/\Z) \ ,&\\
}\]
from which we get an exact sequence $0\rightarrow \br(V)\rightarrow \br(U)\rightarrow H^1(D_{sm},\Q/\Z)$.
Passing to the limit over $D$ gives the Proposition.\qed

\medskip	

Let $K$ be a field. For each discrete valuation ring (DVR) $R$ with quotient field $K$ and residue field $\kappa _R$, there is a natural  exact sequence \cite[III, Prop.\ 2.1]{G10} :
\[0\rightarrow \br(R)\rightarrow \br(K)\qfl{\rho _R}H^1_{\textrm{\'et}}(\kappa _R, \Q/\Z)\ .\]
The group $\br_{\textrm{nr}}(K)$ is defined as the intersection of the subgroups $\Ker \rho _R$, where $R$ runs through all  DVR  with quotient field $K$. 

Now consider the exact sequence of Proposition \ref{sexB}. The group $\varinjlim\limits_U \br(U)$ can be identified with the Brauer group $\br(\C(V))$. For $D\in \mathcal{D}$,  the group $H^1(D_{sm},\Q/\Z)$ embeds into the cohomology group 
$H^1_{\textrm{\'et}}(\C(D), \Q/\Z)$, and  the composition $\br(\C(V))\rightarrow H^1_{\textrm{\'et}}(\C(D), \Q/\Z)$ is equal to the homomorphism $\rho^{} _{\O_{V,D}}$ associated to the DVR $\O_{V,D}$. Thus we have $\br_{\textrm{nr}}(\C(V))\subset \br(V)$. But if $R$ is any DVR with quotient field $\C(V)$, the inclusion $\Spec \C(V)\mono V$ factors as $\Spec \C(V)\mono \Spec R\rightarrow V$ by the valuative criterion of properness, hence $\br(V)$ is contained in the image of $ \br(R)$ in $\br(K)$, that is, in $\Ker \rho _R$. Thus we have $\br(V)=\br_{\textrm{nr}}(\C(V))$ as claimed.

\smallskip	
The big advantage of  working with   $\br_{\textrm{nr}}(K)$  is that to compute it, we do not need to find a smooth projective model of the function field $K$. This was used first by Saltman to give his celebrated counter-example to the Noether problem \cite{Sa} : there exists a finite group $G$ and a linear representation $V$ of $G$ such that the variety $V/G$ is not rational. In such a situation Bogomolov has given a very explicit formula for $\br_{\textrm{nr}}(\C(V/G))$ in terms of the Schur multiplier of $G$ \cite{Bogo}.

\smallskip	
The idea of using the unramified Brauer group to prove non-rationality results has been extended to higher \emph{unramified cohomology groups}, starting with the paper \cite{CTO}. We refer to  \cite{CT} for a survey about these more general invariants.

\medskip	
\section{The Chow group of $0$-cycles}\label{chow}

In this section we discuss another property of (stably) rational varieties, namely the fact that their Chow group $CH_0$ parametrizing $0$-cycles is \emph{universally trivial}.
While the idea goes back to the end of the 70's (see \cite{Bl}), its use for rationality questions is recent \cite{V4}. 

This property implies that $H^3(X,\Z)$ is torsion free, but not conversely. Moreover it behaves well under deformation, even if we accept mild singularities (Proposition \ref{Vdef} below). 
 
In this section we will need to work over non-algebraically closed fields (of characteristic $0$). We use the language of schemes.

Let $X$ be a smooth algebraic variety over a field $k$, of dimension $n$. Recall that the Chow group $CH^p(X)$ is the group of codimension $p$ cycles on $X$  modulo linear equivalence. More precisely,  let us denote by $\Sigma ^p(X)$ the set of codimension $p$ closed subvarieties of $X$. Then $CH^p(X)$ is defined by the exact sequence
\begin{equation}\label{defchow}
\bigoplus_{W\in \Sigma ^{p-1}(X)}k(W)^*\longrightarrow \Z^{(\Sigma ^{p}(X))}\longrightarrow CH^p(X)\rightarrow 0\ ,
\end{equation}
where  the first arrow associates to $f\in k(W)^*$ its divisor \cite[1.3]{Fu}.

We will be particularly interested in the group $CH_0(X):=CH^n(X)$ of $0$-cycles. Associating to  a $0$-cycle $\sum n_i[p_i]$ ($n_i\in\Z, p_i\in X$) the number $\sum n_i[k(p_i):k]$ defines a  homomorphism $\deg: \allowbreak CH_0(X)\rightarrow \Z$. We denote its kernel by $CH_0(X)_0$.

\begin{prop}\label{CH0}
Let $X$ be a smooth  complex projective variety, of dimension $n$, and let $\Delta _X\subset X\times X$ be the diagonal.  The following conditions are equivalent :

{\rm (i)} For every extension $\C\rightarrow K$, $CH_0(X_K)_0=0$;

{\rm (ii)} $CH_0(X_{\C(X)})_0=0$;

{\rm (iii)} There exists a point $x\in X$ and a nonempty Zariski open subset $U\subset X$ such that $\Delta _X-X\times \{x\}$ restricts to $0$ in $CH(U\times X)$;

{\rm (iv)} there exists a point $x\in X$, a smooth projective variety $T$ of dimension $<n$ (not necessarily connected), a generically injective map $i:T\rightarrow X$, and 
a cycle class $\alpha \in CH(T\times X)$ such that 
\begin{equation}\label{diag}
\Delta _X-X\times \{x\} =(i\times 1)_*\alpha \quad\mbox{in }\ CH(X\times X) \ .
\end{equation}

\smallskip	
When these properties hold, we say that $X$ is $CH_0$\emph{-trivial}.
\end{prop} 
\pr The implication  (i) $\Rightarrow$ (ii) is clear.

\smallskip	
\noindent  (ii) $\Rightarrow$ (iii) : Let $\eta $ be the generic point of $X$. The point $(\eta ,\eta )$ of $\{\eta\} \times X =X_{\C(X)}$ is  rational (over $\C(X)$), hence is linearly equivalent to $(\eta ,x)$ for any closed point $x\in X$. The class $\Delta _X-X\times \{x\}$ restricts to $(\eta ,\eta )-(\eta ,x)$ in $CH_0(\eta \times X)$, hence to $0$. We want to show that this implies (iii).
 
 An element of $\Sigma ^p(\eta \times X)$ extends to an element of  $\Sigma ^p(U\times X)$ for some Zariski open subset $U$ of $X$; in other words, the natural map $\varinjlim\limits_U  \ \Sigma ^p(U\times X)\rightarrow \Sigma ^p(\eta \times X)$ is an isomorphism.
Thus writing down the  exact sequence (\ref{defchow}) for $U\times X$ and passing to the direct limit over $U$ we get a commutative diagram of exact sequences
\[\xymatrix@M=8pt{\varinjlim\limits_U \bigoplus\limits_{W\in \Sigma ^{p-1}(U\times X)}k(W)^*\ar[r]\ar@<4ex>[d] &\varinjlim\limits_U\, \Z_{}^{(\Sigma ^{p}(U\times X))}\ar[r]\ar@<4ex>[d]&\varinjlim\limits_U CH^p(U\times X) \ar[r]\ar[d]&0\\
\bigoplus\limits_{W\in \Sigma ^{p-1}(\eta \times X)}k(W)^*\ar[r]&\Z_{}^{(\Sigma ^{p}(\eta \times X))}\ar[r] &CH^p(\eta \times X) \ar[r]&0
}\]
where the  first two vertical arrows are isomorphisms; therefore the third one is also an isomorphism. We conclude that the  class $\Delta- X\times \{x\}$ is zero in $CH^n(U\times X)$ for some $U$.

\medskip	
\noindent (iii) $\Rightarrow$ (iv) : Put $T':=X\smallsetminus U$. The localization exact sequence \cite[Prop.\ 1.8]{Fu}
\[CH(T'\times X)\rightarrow CH(X\times X)\rightarrow CH(U\times X)\rightarrow 0\] 
implies that  $\Delta- X\times \{x\}$ comes from the class in $CH(T'\times X)$ of a cycle  $\sum n_i Z'_i$. For each $i$, let $T'_i$ be the image of $Z_i$ in $T'$, and let $T_i$ be a desingularization of $T'_i$. Since $Z'_i$ is not contained in the singular locus $\Sing(T_i)\times X$, it is the pushforward of an irreducible subvariety $Z_i\subset T_i\times X$. Putting $T=\coprod T_i$ and $\alpha =\sum n_i[Z_i]$ does the job.

 \medskip	
\noindent (iv) $\Rightarrow$ (i) : Assume that  (\ref{diag})  holds; then it holds
  in $CH(X_K\times X_K)$ for any extension $K$ of $\C$, so it suffices to prove  $  CH_0(X)_0=0$.

Denote by $p$ and $q$ the two projections from $X\times X$ to $X$, and put $n:=\dim(X)$. Any class $\delta \in CH^n(X\times X)$ induces a homomorphism $\delta _*:CH_0(X)\rightarrow CH_0(X)$, defined by $\delta _*(z)= q_*(\delta \cdot p^*z)$. Let us consider the classes which appear in (\ref{diag}). The diagonal induces the identity of $CH_0(X)$; the class of $X\times \{x\} $ maps $z\in CH_0(X)$ to $\deg(z)\,[x]$, hence is $0$ on $CH_0(X)_0$. 

Now consider $\delta :=(i\times 1)_*\alpha$. Let $p',q'$ be the projections from $T\times X$ to $T$ and $X$. Then, for $z\in CH_0(X)$,
\[\delta _*z=q_*((i\times 1)_*\alpha \cdot p^*z)=q'_*(\alpha \cdot p'^*i^*z)\ .\]
Since $\dim T<\dim X$, $i^*z$ is zero, hence also $\delta _*z$.
We conclude from (\ref{diag}) that  $CH_0(X)_0=0$. \qed

\medskip	
\noindent \emph{Example}$.-$ The group $CH_0(X)$ is a birational invariant \cite[ex.\ 16.1.11]{Fu}, thus the above properties depend only on the birational equivalence class of $X$. In particular  a rational variety  is $CH_0$-trivial. 
More generally, since $CH_0(X\times \P^n)\cong CH_0(X)$ for any variety $X$, a stably rational variety is $CH_0$-trivial. 

\medskip	
Despite its technical aspect, Proposition \ref{CH0}
 has remarkable consequences, which have been worked out by Bloch and Srinivas \cite{BS} :   

\begin{prop}\label{BS}
Suppose $X$ is $CH_0$-trivial.

$1)$ $H^0(X,\Omega ^r_X)=0$ for all $r>0$.

$2)$ The  group $H^3(X,\Z)$ is torsion free.
\end{prop}
\pr  The proof is very similar to
 that of the implication (iii)$\ \Rightarrow\ $(i) in the previous Proposition; we use the same notation. Again a class $\delta $ in $CH^n(X\times X)$ induces a homomorphism $\delta^*: H^r(X,\Z)\rightarrow H^r(X,\Z) $,  defined by $\delta ^*z:=p_*(\delta \cdot q^*z)$. The diagonal induces the identity, the class $[X\times \{p\} ]$ gives $0$ for $r>0$, and the class $(i\times 1)_*\alpha $ gives the homomorphism $z\mapsto i_*p'_*(\alpha \cdot q'^*z)$. Thus  formula (\ref{diag}) gives for $r>0$ a commutative diagram
 \begin{equation}\label{tri}
\xymatrix{& H^*(T,\Z)\ar[dr]^{i_*}&\\ H^r(X,\Z)\ar[ru]\ar[rr]^{\mathrm{Id}} &&H^r(X,\Z)\ .
 }\end{equation}
  On each component $T_k$ of $T$ the homomorphism $i_*: H^*(T_k,\C)\rightarrow H^*(X,\C)$ is a morphism of Hodge structures of bidegree $(c,c)$, with $c:=\dim X-\dim T_k>0$. Therefore its image intersects trivially the subspace $H^{r,0}$ of $H^r(X,\C)$. 
Since $i_*$ is surjective by (\ref{tri}), we get $H^{r,0}=0$.

Now we take $r=3$ in (\ref{tri}). The only  part of $H^*(T,\Z)$ with a nontrivial contribution in (\ref{tri}) is $H^1(T,\Z)$, which is torsion free. Any torsion element in $H^3(X,\Z)$ goes to $0$ in $H^1(T,\Z)$, hence is zero.\qed

 \medskip	
 Observe that  in the proof we use only formula (\ref{diag}) in $H^*(X\times X)$ and not in the Chow group. The relation between these two properties is discussed in Voisin's papers \cite{V3,V4,V5}.

\medskip	
As the Clemens-Griffiths criterion, the  triviality of $CH_0(X)$ behaves well under deformation (compare with Lemma \ref{defjac}) :

\begin{prop}\cite{V4}\label{Vdef}
 \ Let $\pi :X\rightarrow B$ be a flat, proper family over a smooth variety $B$, with $\dim(X)\geq 3$. Let $\mathrm{o}\in B$; assume that :

$\bullet$ The general fiber $X_b$ is smooth;

 $\bullet$ $X_{\mathrm{o}}$ has only ordinary double points, and its   
  desingularization $\tilde{X}_{\mathrm{o}} $   is not $CH_0$-trivial. 
  
 Then    ${X} _b$ is not $CH_0$-trivial for  a very general point $b$ of $B$.\end{prop}
 Recall that  \`{}very general' means \`{}outside a countable union of strict subvarieties of $B$' (\ref{K}).
 
We refer to \cite{V4} for the proof. The idea is that there cannot exist a decomposition (\ref{diag}) of Proposition \ref{CH0} for $b$ general in $B$, because it would extend to an analogous decomposition over $X$, then specialize to $X_{\mathrm{o}}$, and finally extend to $\tilde{X}_{\mathrm{o}} $. One concludes by observing that the locus of points $b\in B$ such that $X_b$ is smooth and $CH_0$-trivial is a countable union of  subvarieties.
 \begin{cor}
The double cover of $\P^3$ branched along a very general quartic surface is not stably rational.
\end{cor}
\pr Consider the pencil of quartic surfaces in $\P^3$ spanned by a smooth quartic and a quartic symmetroid, and  the family of double covers of $\P^3$ branched along the members of this pencil. By Proposition \ref{BS}.2), the Artin-Mumford threefold is not $CH_0$-trivial.  Applying the Proposition we conclude that 
  a very general quartic double solid is not $CH_0$-trivial, hence not stably rational.
\qed

\smallskip	
More generally, Voisin shows that the desingularization of a very general quartic double solid with at most seven nodes is not stably rational. 

\medskip	
Voisin's idea has given rise to a number of new results. Colliot-Th\'el\`ene and Pirutka have extended Proposition \ref{Vdef} to the case where the singular fiber $X_{\mathrm{o}}$ has (sufficiently nice) non-isolated singularities, and applied this to prove that a very general quartic threefold is not stably rational \cite{ACT}. Using their extension and Proposition \ref{2P} I have shown that the double cover of $\P^4$ or $\P^5$ branched along a very general quartic hypersurface is not stably rational \cite{B2Q}, and that a very general sextic double solid is not stably rational \cite{B2S}. As already mentioned, combining the methods of Koll\'ar and Voisin, Totaro has proved that a very general hypersurface of degree $d$ and dimension $n$ is not stably rational for $d\geq 2\lceil \frac{n+2}{3}\rceil  $ \cite{T}. Finally, Hassett, Kresch and Tschinkel have shown  that a conic bundle (see \ref{schott}) with discriminant   a very general plane curve of degree $\geq 6$   is not stably rational \cite{HKT}.

\medskip	
We do not know whether there exist smooth quartic double solids which are $CH_0$-trivial. In contrast, Voisin has constructed families of smooth cubic threefolds wich are $CH_0$-trivial \cite{V5} -- we do not know what happens for a general cubic threefold.

\end{document}